\documentclass[11pt]{amsart}

\usepackage[article]{robertpreamble}
\usepackage[margin=1in]{geometry}
\usepackage{booktabs}

\newcommand{\HP}[1][k]{\mathcal H^{#1}_P}
\newcommand{\dP}{d_P}
\newcommand{\Hor}{\mathcal H}
\newcommand{\Ver}{\mathcal V}
\newcommand{\Adm}{\mathcal A}
\newcommand{\Fam}{\mathcal T}
\newcommand{\restr}{\mathbin{\vrule height 1.6ex depth 0pt width 0.13ex\vrule height 0.13ex depth 0pt width 1.1ex}}
\newcommand{\BMOP}{\operatorname{BMO}_P}
\newcommand{\Dt}{|D_t|^{1/2}}
\newcommand{\releasecommit}{aef12f7c3f3e4aa5b54cbad57f9cd8dd7148226e}
\DeclareMathOperator{\LipP}{Lip_P}
\definecolor{setcol}{RGB}{54,79,171}
\definecolor{modelcol}{RGB}{245,132,39}
\colorlet{zoomcol}{modelcol}

\begin{document}

\title{Parabolic rectifiable Reifenberg for measures}
\author{Robert Koirala}
\address{Department of Mathematics, University of California San Diego}
\email{rkoirala@ucsd.edu}
\date{\today}

\begin{abstract}
    We announce a structure theory for measures on $\R^n\times\R$ equipped with the parabolic distance. Our results include a parabolic rectifiable Reifenberg theorem for arbitrary Radon measures, qualitative and quantitative characterizations of regular parabolic rectifiability by Jones square functions, and sharp results for $L^p$ beta numbers. The results hold in every parabolic dimension and extend the theories of Edelen--Naber--Valtorta, Azzam--Tolsa and David--Semmes to space-time. The stated theorems and corollaries, together with their proofs, are formalized in Lean~4 with Mathlib.
\end{abstract}

\maketitle

\vspace{-4mm}
\section{Introduction}

In the work announced here \cite{KoiralaPrep}, we study the structure of measures on parabolic space-time $\R^n\times\R$ that look flat at most locations, scales and resolutions, as measured by parabolic Jones beta numbers \cite{Jones1990}. The beta number $\beta_\mu(z,r)$ is a normalized $L^2$-distance from $\mu$ to flat models in the parabolic ball $P(z,r)$, and these models are invariant under the parabolic scaling $(x,t)\mapsto(rx,r^2t)$ (see \eqref{eq:papa-jones}). Our main results are the following.

\begin{enumerate}
    \item \emph{Parabolic rectifiable Reifenberg theorem} (Theorem~\ref{thm:reifenberg}). Under a weak bound on the square function $S_\mu$ of the beta numbers, a locally Radon measure $\mu$ decomposes into a part carried by graphs over horizontal flat models, a part carried by \textit{regular} graphs over vertical flat models, and a remainder of controlled mass. The two graph parts satisfy uniform content and packing estimates, and the graphs have controlled regularity.

    \item \emph{Qualitative characterization} (Theorem~\ref{thm:qualitative}). A measure $\mu$ is \textit{regularly parabolic rectifiable} if and only if $S_\mu$ is finite and the upper $k$-density of $\mu$ is positive and finite almost everywhere. Under these density conditions, $\mu$ also decomposes uniquely into horizontal, vertical and purely unrectifiable parts (Corollary~\ref{cor:canonical}).

    \item \emph{Quantitative characterization} (Theorem~\ref{thm:quantitative}). For Ahlfors regular
    measures, a Carleson square-function bound is equivalent to a uniform covering of points and scales by horizontal and regular vertical graphs (Section~\ref{sec:quantitative}).

    \item \emph{Different exponents} (Corollary~\ref{cor:table}). We determine for which exponents $p$ the square functions of $L^p$ beta numbers are necessary or sufficient for regular rectifiability. As in Tolsa's Euclidean work \cite{Tolsa2019beta}, only $p=2$ gives a characterization.
\end{enumerate}

\begin{remark}
\begin{enumerate}
    \item Theorem~\ref{thm:reifenberg} significantly extends the Euclidean Reifenberg theory of Naber--Valtorta and Edelen--Naber--Valtorta \cite{ENV2025,NaberValtorta2017} to space-time. The main obstacle is conceptual: it comes from the anisotropy of time and from the regularity required in the time direction. In our joint work with Hallgren and Ma on parabolic nodal and singular sets \cite{HallgrenKoiralaMa25}, the models were vertical and we overcame this obstacle using quantitative uniqueness of the approximating models; other recent works in the parabolic setting also rely on some form of quantitative uniqueness \cite{Bamler20a,CheegerHaslhoferNaber13,CheegerHaslhoferNaber14,ColdingMinicozzi16,FangLi25a,FangLi25b,FuWangWuZhang26,Gianniotis19,LW26}. Techniques that work without quantitative uniqueness have been missing, and our work removes the obstacle in that case.
    
    \item The qualitative characterization is a parabolic counterpart of the theorems of Azzam--Tolsa and Tolsa \cite{AzzamTolsa2015,Tolsa2015}. It extends the codimension-one graph characterization of Bortz, Egert, Ferris and Saari \cite{BEFS26} to arbitrary measures and to both types of models; in codimension one only vertical models occur. The quantitative characterization is a parabolic analogue of the results of Jones and David--Semmes \cite{DavidSemmes1991,DavidSemmes1993,Jones1990}.
    
    \item Regular parabolic rectifiability uses countably many horizontal graphs and regular vertical graphs. The latter are parabolic Lipschitz graphs whose half-time derivative lies in parabolic $\operatorname{BMO}$; see Section~\ref{sec:decomposition} for the precise definitions. Regularity is strictly stronger than the parabolic Lipschitz property (Example~\ref{ex:lacunary} and \cite{Mat22}). The role of regular graphs in parabolic rectifiability was suggested by work on caloric removability, parabolic uniform rectifiability and nodal sets \cite{BEFS26,BHHLN23JGA,HallgrenKoiralaMa25,HLN03,MPT22}.
    
    \item The hypothesis of Theorem~\ref{thm:reifenberg} is much weaker than those of the characterization theorems and makes no assumption on the measure, yet it gives a decomposition of the same kind, at the cost of a controlled remainder. Without Ahlfors regularity such a decomposition has no converse, since beta numbers measure mass as well as flatness (Remark~\ref{rem:converse}). Weak square-function hypotheses of this kind arise from monotone quantities in geometric flows, where we expect applications as in the elliptic theory \cite{ENV2025,NaberValtorta2017}. For instance, Theorem~\ref{thm:reifenberg} upgrades the parabolic rectifiability of all the quantitative strata in \cite{HallgrenKoiralaMa25} to regular parabolic rectifiability; there we proved regular parabolic rectifiability of the nodal and singular sets. Such applications to geometric flows will be discussed elsewhere.
\end{enumerate}
\end{remark}

All theorems and corollaries stated in Sections~\ref{sec:decomposition}--\ref{sec:betap}, together with their proofs, are formalized in Lean~4 with Mathlib and checked against a Mathlib-only specification (Section~\ref{sec:lean}). A detailed exposition, further results and examples, and a discussion of sharpness will appear in \cite{KoiralaPrep}.

\subsection{Historical background}
The theory of rectifiability in Euclidean space is well developed, beginning with the work of Besicovitch, Federer, Marstrand, Mattila and Preiss \cite{Besicovitch1928,Federer1947,Marstrand1964,Mattila1975,Preiss1987}; see the surveys \cite{MatSurvey,Naber2018,Tolsa26ICM}. In 1990, Jones introduced beta numbers in his solution of the analytic traveling salesman problem: they control the length of a rectifiable curve passing through a set \cite{Jones1990,Okikiolu1992,Schul2007}. David and Semmes then characterized uniform rectifiability of sets by Carleson estimates on beta numbers, by corona decompositions and by the boundedness of singular integrals \cite{DavidSemmes1991,DavidSemmes1993}. For measures, the work of Pajot and Badger--Schul \cite{BadgerSchul2016,BadgerSchul2017,Pajot1997} preceded the square-function characterization of rectifiability by Azzam--Tolsa and Tolsa \cite{AzzamTolsa2015,Tolsa2015,Tolsa2019corr}, and the special role of the exponent two was established in \cite{Tolsa2019beta}.

In parallel, Reifenberg proved his topological disk theorem to solve the Plateau problem for
surfaces of arbitrary dimension \cite{Reifenberg1960}: a closed set that is close to a plane at every location and scale is a topological disk. Toro and David--Toro then obtained bi-Lipschitz parametrizations from square-summable flatness, also for sets with holes \cite{DavidToro1999,DavidToro2012,Toro1995}. Naber and Valtorta proved rectifiable and discrete versions of Reifenberg's theorem \cite{Miskiewicz2018,NaberValtorta2017}, and Edelen, Naber and Valtorta treated arbitrary measures under a weak square-function bound \cite{ENV2019,ENV2025}; see also \cite{Naber2018}. These tools have had many applications to the structure of singular sets in nonlinear elliptic problems:
harmonic maps, minimal surfaces and the Yang--Mills equations \cite{Simon1993,Simon1995,CheegerNaber13,NaberValtorta2017,naber-valtorta-YM,NaberValtorta20,naber-valtorta-2024energy,DeLellisMarcheseSpadaroValtorta18},
critical and nodal sets of elliptic equations
\cite{CheegerNaberValtorta15,NaberValtorta17b,HuangJiang23}, free boundaries
\cite{EdelenEngelstein19}, and spaces with Ricci curvature bounds \cite{CheegerNaber12,CheegerNaber15,JiangNaber21,CheegerJiangNaber21}. There are related developments in the parabolic setting, including nodal and singular sets, mean curvature flow, harmonic map flow and Ricci flow \cite{CheegerHaslhoferNaber13,CheegerHaslhoferNaber14,ColdingMinicozzi16,Gianniotis19,Bamler20a,FangLi25a,FangLi25b,FuWangWuZhang26,HallgrenKoiralaMa25,HuangJiang24,LW26}. These works, however, rely on some form of quantitative uniqueness of tangent flows to obtain sharp estimates. Parabolic rectifiable Reifenberg tools for settings without quantitative uniqueness were not available before our work.

In contrast, the theory of rectifiability in the parabolic setting is still at an early stage, as Mattila notes \cite{MatSurvey}. The main reason is the anisotropy of time, which is a conceptual difficulty and not merely a technical one. Several analogues of the Euclidean theory fail: a Preiss-type density theorem fails for the max metric $\dP$ in \eqref{eq:max-metric} \cite{MMP26}; a singular-integral characterization in the spirit of David--Semmes fails \cite{HJ26}; a parabolic rectifiable set can meet every regular parabolic Lipschitz graph in a null set \cite[Example~8.2]{Mat22}; and on general parabolic Lipschitz graphs caloric measure may be singular \cite{KW80,LS88}. The last phenomenon led to the notion of regular parabolic Lipschitz graphs in the layer-potential theory for the heat equation \cite{HL96,Hof97,LM95}. Since then there has been progress in codimension one: parabolic uniform rectifiability has been developed through big pieces of regular parabolic Lipschitz graphs and corona decompositions \cite{BHHLN22AIF,BHHLN22AFM,BHHLN23JGA,BHS26,HLN03,HLN04,NS17}, with connections to the Dirichlet problem and to estimates for caloric functions \cite{BHHW26,BHHLN23CME,BHMN25graph,BHMN25pre}. The work of Mateu, Prat and Tolsa on removable singularities for caloric functions \cite{MPT22}, the work on parabolic uniform rectifiability \cite{HLN03,BHHLN23JGA,BEFS26}, and our joint work with Hallgren and Ma on nodal and singular sets \cite{HallgrenKoiralaMa25} further suggested that integral beta-number conditions should be equivalent to regular parabolic rectifiability. Our work establishes this for arbitrary measures in all dimensions and for both types of models, deriving and characterizing regularity from beta numbers without assuming uniqueness of tangent models.

\subsection{Main difficulties}\label{sec:difficulties}
The main conceptual difficulties come from the interaction of the parabolic geometry, given by the metric
\[
 \dP((x,t),(y,s))=\max\{|x-y|,|t-s|^{1/2}\}
\]
with multiscale approximation and with regularity in time. Further difficulties and their resolution will be discussed in detail in \cite{KoiralaPrep}.

\begin{enumerate}[label=(D\arabic*)]
\item \emph{Anisotropy in time.} Space scales by $r$ and time by $r^2$, so small Euclidean errors in time become large parabolic errors at small scales and can change the dimension. For instance, the line $\{(x,ax)\}$ in $\R\times\R$ with small $a>0$ looks horizontal at coarse scales, but it has parabolic dimension two and vertical tangents (Example~\ref{ex:slanted}). Spatial coherence of the approximating horizontal planes at coarse scales does not control this effect; one also needs quadratic control of the time differences. Similarly, vertical charts need $1/2$-H\"older regularity in time to be parabolic Lipschitz.

\begin{figure}[!htb]
\centering\begin{tikzpicture}[>=Stealth,font=\scriptsize,scale=0.85]
 \begin{scope}
  \draw[gray!60] (-1,-1) rectangle (1,1);
  \draw[->,gray] (-1.1,0)--(1.2,0) node[right] {$x$};
  \draw[->,gray] (0,-1.1)--(0,1.2) node[above] {$t$};
  \draw[very thick,setcol] (-0.85,0.08)--(-0.2,0.08)--(0.2,-0.08)--(0.85,-0.08);
  \draw[thick,zoomcol] (-0.15,-0.0225) rectangle (0.15,0.0225);
  \node at (0,-1.4) {$P(0,1)$};
 \end{scope}
 \draw[zoomcol] (0.15,0.0225)--(2.0,1); \draw[zoomcol] (0.15,-0.0225)--(2.0,-1);
 \begin{scope}[xshift=3cm]
  \draw[zoomcol] (-1,-1) rectangle (1,1);
  \draw[->,gray] (-1.1,0)--(1.2,0) node[right] {$u$};
  \draw[->,gray] (0,-1.1)--(0,1.2) node[above] {$s$};
  \draw[dashed,modelcol] (0,-1)--(0,1);
  \draw[very thick,setcol] (-0.375,1)--(0.375,-1);
  \node at (0,-1.4) {$P(0,\rho)$, rescaled};
 \end{scope}
\end{tikzpicture}
\caption{(D1) A ramp with a small time tilt is close to a horizontal line in $P(0,1)$. In the small ball $P(0,\rho)$, after parabolic rescaling, it is close to
the vertical line (dashed), so it has parabolic dimension two.}\label{fig:D1}
\end{figure}
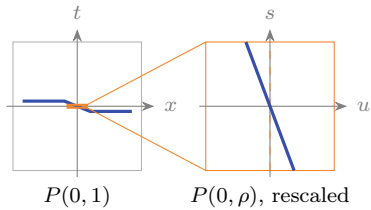

\item \emph{Accumulation of plane errors.} Best approximating planes can be non-unique and can rotate across scales and locations, unlike the vertical planes in \cite{HallgrenKoiralaMa25}; the same situation arises in nonlinear geometric flows. Without summable control, repeated small errors may only give an $\alpha$-H\"older chart with $\alpha<1$, which yields the weaker dimension bound $k/\alpha$ and no sharp $k$-dimensional measure estimate.

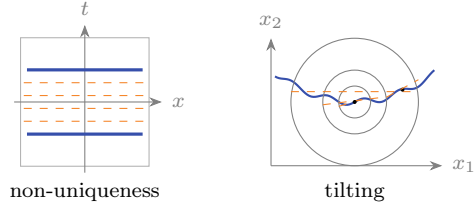
\begin{figure}[!htb]
\centering\begin{tikzpicture}[>=Stealth,font=\scriptsize,scale=0.85]
 \begin{scope}
  \draw[gray!60] (-1,-1) rectangle (1,1);
  \draw[->,gray] (-1.1,0)--(1.2,0) node[right] {$x$};
  \draw[->,gray] (0,-1.1)--(0,1.2) node[above] {$t$};
  \draw[very thick,setcol] (-0.9,0.5)--(0.9,0.5);
  \draw[very thick,setcol] (-0.9,-0.5)--(0.9,-0.5);
  \foreach \y in {-0.3,-0.1,0.1,0.3} \draw[dashed,modelcol] (-0.95,\y)--(0.95,\y);
  \node at (0,-1.4) {non-uniqueness};
 \end{scope}
 \begin{scope}[xshift=4.2cm,scale=0.62]
  \draw[->,gray] (-2.1,-1.6)--(2.2,-1.6) node[right] {$x_1$};
  \draw[->,gray] (-2.1,-1.6)--(-2.1,1.6) node[above] {$x_2$};
  \draw[thick,setcol,domain=-2:2,samples=80] plot (\x,{0.18*\x*\x+0.08*sin(400*\x)});
  \draw[gray] (0,0) circle (1.6); \draw[gray] (0,0) circle (0.8); \draw[gray] (0,0) circle (0.4);
  \draw[dashed,modelcol] (-1.6,0.26)--(1.6,0.26);
  \draw[dashed,modelcol] (-0.8,-0.07)--(0.8,0.13);
  \draw[dashed,modelcol] (-0.4,-0.11)--(0.4,0.13);
  \fill (0,0) circle (0.05);
  \fill (1.2,0.31) circle (0.05);
  \draw[dashed,modelcol] (0.8,0.12)--(1.6,0.56);
  \node at (0,-2.25) {tilting};
 \end{scope}
\end{tikzpicture}
\caption{(D2) Left: For a measure split equally between the time slices $t=\pm h$, every $\tau\in[-h,h]$ gives a best fit
(dashed). Right, in a spatial slice: the best spatial plane (dashed) rotates as the scale decreases,
and differs at nearby points.}\label{fig:D2}
\end{figure}

\item \emph{Change of model type.} At a fixed $k$, horizontal models have $k$ spatial
 directions and vertical models have $k-2$, so the two types cannot be joined by spatial
 rotations. Repeated changes of type must preserve the content and packing estimates.

\begin{figure}[!htb]
\centering\begin{tikzpicture}[>=Stealth,font=\scriptsize,scale=1.05]
 \begin{scope}[x={(1cm,0cm)},y={(0.45cm,0.3cm)},z={(0cm,1cm)}]
  \draw[gray!60] (-1,-1,-1)--(1,-1,-1)--(1,1,-1)--(-1,1,-1)--cycle;
  \foreach \x/\y in {-1/-1,1/-1,1/1,-1/1} \draw[gray!60] (\x,\y,-1)--(\x,\y,1);
  \draw[very thick,setcol] (0,0,-0.2)--(0,0,0);
  \fill[blue!12,even odd rule] (-1,-1,0)--(1,-1,0)--(1,1,0)--(-1,1,0)--cycle
   plot[domain=0:360,samples=72,smooth] ({0.22*cos(\x)},{0.22*sin(\x)},0);
  \draw[setcol] (-1,-1,0)--(1,-1,0)--(1,1,0)--(-1,1,0)--cycle;
  \draw[setcol] plot[domain=0:360,samples=72,smooth] ({0.22*cos(\x)},{0.22*sin(\x)},0);
  \draw[very thick,setcol] (0,0,0)--(0,0,0.2);
  \draw[gray!60] (-1,-1,1)--(1,-1,1)--(1,1,1)--(-1,1,1)--cycle;
  \draw[thick,zoomcol] (-0.12,0,-0.0144)--(0.12,0,-0.0144)--(0.12,0,0.0144)--(-0.12,0,0.0144)--cycle;
  \draw[zoomcol] (0.12,0,0.0144)--(2.6,0,1); \draw[zoomcol] (0.12,0,-0.0144)--(2.6,0,-1);
  \draw[->,gray] (1,-1,0)--(1.5,-1,0) node[right] {$x_1$};
  \draw[->,gray] (-1,1,0)--(-1,1.7,0) node[above] {$x_2$};
  \draw[->,gray] (0,0,1)--(0,0,1.3) node[above] {$t$};
  \draw[gray,densely dotted] (0,0,0.2)--(0,0,1);
  \node at (0,0,-1.35) {$P(0,1)$};
 \end{scope}
 \begin{scope}[xshift=3.6cm]
  \draw[zoomcol] (-1,-1) rectangle (1,1);
  \draw[->,gray] (-1.1,0)--(1.25,0) node[right] {$x_1/\rho$};
  \draw[->,gray] (0,-1.1)--(0,1.25) node[above] {$t/\rho^2$};
  \draw[very thick,setcol] (0,-1)--(0,1);
  \draw[dashed,modelcol] (0.07,-1)--(0.07,1);
  \node at (0,-1.35) {$P(0,\rho)$, rescaled};
 \end{scope}
\end{tikzpicture}
\caption{(D3) A change of type at a fixed dimension $k=2$ in $\R^2\times\R$: a horizontal
plane with a hole, and a vertical needle $\{0\}\times[-h,h]$ through the hole. In $P(0,1)$ the best
model is the horizontal plane. The small ball $P(0,\rho)$, with $\rho$ smaller than the hole and
$\rho^2<h$, meets only the needle, and there the best model is the vertical line $\{0\}\times\R$
(dashed). Both models have parabolic dimension two.}\label{fig:D3}
\end{figure}
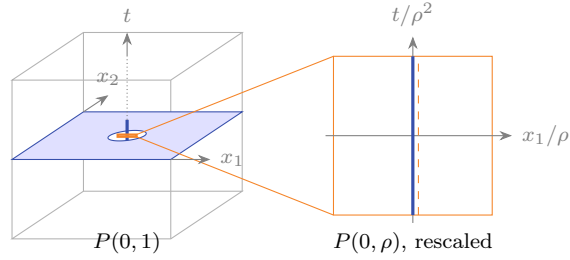

\item \emph{Regularity of the graph covers.} Parabolic Lipschitz control does not imply a
 $\operatorname{BMO}$ half-time derivative (Example~\ref{ex:lacunary}). The construction
 initially produces partial graphs, and extending them to full regular graphs must preserve
 the affine-error estimates that this nonlocal regularity in time requires.

\begin{figure}[!htb]
\centering\includegraphics[width=0.75\textwidth]{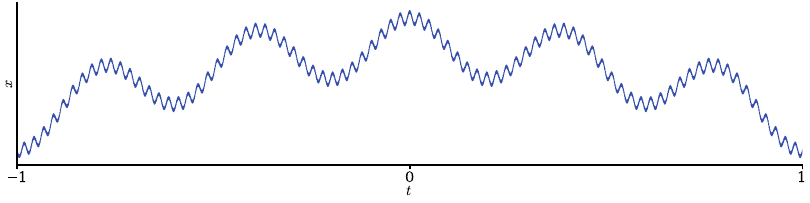}
\caption{(D4) The graph of $f(t)=\frac3{16}\sum_{j\ge0}(j+1)^{-1/2}16^{-j/2}\cos(16^jt)$ drawn with the terms $j\le6$. Here $f$ is parabolic Lipschitz but not regular parabolic Lipschitz.}\label{fig:D4}
\end{figure}
\end{enumerate}

\subsection{Key ideas}\label{sec:ideas}
We develop neck regions and neck decompositions
\cite{CheegerJiangNaber21,ENV2025,JiangNaber21,NaberValtorta2017} adapted to the parabolic
geometry, inspired by \cite{HallgrenKoiralaMa25}. Once the difficulties above are addressed, the architecture of the proof of the parabolic rectifiable Reifenberg theorem is inspired by \cite{ENV2025,Naber2018}, and the characterization theorems are inspired by the work of Azzam--Tolsa, David--Semmes and Tolsa \cite{AzzamTolsa2015,DavidSemmes1991,Tolsa2015}. We now describe some of the key ideas. Throughout, $\mu$ is a locally Radon measure on $\R^n\times\R$, and $P(z,r)$ is the $\dP$-ball of radius $r$ centered at $z\in\R^n\times\R$.

\begin{enumerate}[label=(I\arabic*)]
\item \emph{Coherent space-time charts} (D1--D2). Let $\lambda=\mu\restr P(z,r)$, and let $H=V\times\{\tau\}$ and $H'=V'\times\{\tau'\}$ be two best models in $P(z,r)$, both horizontal, at times $\tau$ and $\tau'$. Since $\dP((x,t),H)^2\ge|t-\tau|$,
\[
    \lambda(\R^n\times \R)\,|\tau-\tau'|\le
    \int\bigl(|t-\tau|+|t-\tau'|\bigr)\,d\lambda \lesssim \beta_\mu(z,r)^2 r^{k+2}.
\]
Thus, if $\lambda(\R^n\times \R)\gtrsim r^k$ and $\beta_\mu(z,r)\le\varepsilon$, then $|\tau-\tau'|\lesssim\varepsilon^2r^2$, a quadratic gain. The collapsed case $\lambda(\R^n\times \R)\ll r^k$ is easier, because the content is small. To handle non-uniqueness, we run a Reifenberg-type argument with multiscale projections and gluing, as in \cite{Naber2018,NaberValtorta2017,Simon1993}, adapted to the parabolic geometry by also tracking the time differences. This also gives $1/2$-H\"older control in the time direction.

\item \emph{Content contraction at a change of type} (D3). The key observation is that a $k$-dimensional horizontal plane and a $k$-dimensional vertical plane intersect in a set contained in a $(k-2)$-dimensional horizontal plane. A quantitative version of this observation controls the content of the regions where the type of the model changes.

\item \emph{Fourier estimates and regular extension} (D4). Let
 $f\in L^2(\R^m\times\R)$ be a function on a vertical base, with variables $(u,t)$, and let
 $\widehat f(\xi,\tau)$ be its Fourier transform, with frequencies $\xi\in\R^m$ in space and
 $\tau\in\R$ in time. Spatially affine approximation satisfies
 \[
  \int|\tau|\,|\widehat f|^2
  \lesssim\iint_0^\infty\alpha_f(p,r)^2\frac{dr}{r}\,dp
  \lesssim\int(|\xi|^2+|\tau|)\,|\widehat f|^2,
 \]
 where $\alpha_f(p,r)^2=r^{-m-4}\inf_{a,b}\int_{P(p,r)}|f-a-b\cdot u|^2$ is the analogue of the beta numbers for functions. These are parabolic versions of Dorronsoro's estimates \cite{Dorronsoro1985}; for Hofmann's estimate
 see \cite{Hof97} and \cite[Theorem~3.8]{Orp21}. Localizing the first inequality, following the criterion in \cite[Proposition~B.1]{HallgrenKoiralaMa25} and \cite{HLN03}, gives $\|\Dt F\|_{\BMOP}\le C(L+\sqrt{B_0})$ from a Lipschitz bound $L$ and a Carleson bound $B_0$ on $\alpha_F$. The partial graphs are decomposed, up to $\HP$-null sets, into countably many pieces with regular extensions and controlled affine errors, as in \cite{HallgrenKoiralaMa25}. After localization, the second inequality gives the converse direction of the characterization on regular graphs.
 
\end{enumerate}

\subsection*{LLM usage}
Large language models were used to write the Lean code and documentation of the formalization under the author's direction, and to assist with the literature search and the drafting of this announcement. The author is responsible for the mathematics, the architecture/sketches of the proof and formalization, the definitions and statements, and this text. Every theorem and corollary in Sections~\ref{sec:decomposition}--\ref{sec:betap}, and Example~\ref{ex:lacunary}, is a translation of a Lean theorem, and all their proofs are machine-checked (Section~\ref{sec:lean}). Examples~\ref{ex:slanted} and~\ref{ex:tilt} and Remark~\ref{rem:converse} are elementary illustrations checked by hand and are not formalized. The description of the ideas in
Section~\ref{sec:difficulties}--\ref{sec:ideas} is informal.

\subsection*{Acknowledgement}
The author thanks Bennett Chow and Ruobing Zhang for their constant support and encouragement. He is very grateful to Yuan Liao and Ziyang Qin for discussion on autoformalization, and especially Yuan Liao for his guidance early in the process,
which sped up the machine verification of the results announced here.

\section{Decomposition for measures}\label{sec:decomposition}

Fix an integer $n\ge1$ and equip $\R^n\times\R$ with the metric
\begin{align}\label{eq:max-metric}
    \dP((x,t),(y,s))=\max\{|x-y|,|t-s|^{1/2}\},
\end{align}
the parabolic scaling $(x,t)\mapsto(rx,r^2t)$ and Hausdorff dimension $n+2$. We write $P(z,r)$ for the
open $\dP$-ball, $\HP$ for the (unnormalized) $k$-dimensional Hausdorff measure of $\dP$ and
$|\cdot|$ for Lebesgue measure on $\R^n\times\R$, and $\mathcal{L}^k$ for Lebesgue measure on Euclidean $\R^k$. All dimensions $k$ are integers with $1\le k\le n+2$. A
measure means a nonnegative Borel measure; it is \emph{locally Radon} if it is finite on compact
sets and inner regular. Its total mass may be infinite.

\subsection{Models, beta numbers and graphs}
The \emph{horizontal} models $\Hor_k$ are the sets $(a+V)\times\{t_0\}$ with $V\subset\R^n$ a linear
subspace, $\dim V=k$; they exist for $1\le k\le n$. The \emph{vertical} models $\Ver_k$ are the sets
$(a+W)\times\R$ with $\dim W=k-2$; they exist for $2\le k\le n+2$, and for $k=n+2$ the only vertical
model is $\R^n\times\R$. Let $\Adm_k=\Hor_k\cup\Ver_k$, and let $\Fam$ denote one of the families
$\Hor_k,\Ver_k,\Adm_k$. For a measure $\mu$, a family
$\Fam$, $z\in\R^n\times\R$ and $r>0$, the parabolic beta number and its square function are
\begin{align}\label{eq:papa-jones}
    \beta_{\mu,\Fam}(z,r)^2=\frac1{r^{k+2}}\inf_{L\in\Fam}\int_{P(z,r)}\dP(w,L)^2\,d\mu(w),\qquad S^R_{\mu,\Fam}(z)=\int_0^R\beta_{\mu,\Fam}(z,r)^2\frac{dr}{r}\in[0,\infty].
\end{align}
The superscript $R$ denotes the largest scale, not a power. We write $\beta_\mu=\beta_{\mu,\Adm_k}$ and $S^R_\mu=S^R_{\mu,\Adm_k}$. The lower and upper densities are
$\Theta^k_*(\mu,z)=\liminf_{r\downarrow0}r^{-k}\mu(P(z,r))$ and
$\Theta^{k,*}(\mu,z)=\limsup_{r\downarrow0}r^{-k}\mu(P(z,r))$.

Let $W\subset\R^n$ be a linear subspace. A map $G:W\times\R\to W^\perp$ is \emph{parabolic Lipschitz}
with constant $L$ if $|G(w,t)-G(w',s)|\le L\,\dP((w,t),(w',s))$. Its \emph{half-time derivative}
$\Dt G$ is the Fourier multiplier $|\tau|^{1/2}$ in the time variable,
\[
 \widehat{\Dt G}(w,\tau)=|\tau|^{1/2}\,\widehat G(w,\tau),
\]
where $\widehat{\ }$ is the Fourier transform in $t$; for a parabolic Lipschitz $G$ it is defined in the
sense of distributions, componentwise and modulo constants. A function $g$ on $W\times\R$ has
\emph{parabolic $\operatorname{BMO}$ seminorm}
\[
 \|g\|_{\BMOP}=\sup_{P}\ \frac1{|P|}\int_P|g-g_P|,\qquad g_P=\frac1{|P|}\int_Pg,
\]
the supremum being over all $\dP$-balls $P$ of $W\times\R$.

\begin{definition}\label{def:graphs}
Let $\varepsilon>0$.
\begin{enumerate}
\item A \emph{horizontal graph of slope at most $\varepsilon$} is a set $\{(a+u+f(u),\tau(u)):u\in D\}$,
 where $V\subset\R^n$ is a $k$-dimensional subspace, $D\subset V$ is an \emph{arbitrary} subset,
 $f:D\to V^\perp$, $\tau:D\to\R$ and
 \[\max\{|f(u)-f(v)|,|\tau(u)-\tau(v)|^{1/2}\}\le\varepsilon|u-v|\qquad(u,v\in D).\]
\item A \emph{regular vertical graph of norm at most $\varepsilon$} is a set
 $\{(a+w+G(w,t),t):(w,t)\in W\times\R\}$, where $W\subset\R^n$ is a $(k-2)$-dimensional subspace and
 $G:W\times\R\to W^\perp$ is parabolic Lipschitz with constant $L$ and
 $L+\|\Dt G\|_{\BMOP}\le\varepsilon$.
\end{enumerate}
\end{definition}

Because $|\tau(u)-\tau(v)|\le\varepsilon^2|u-v|^2$, the time function of a horizontal graph is locally
constant on the interior of its domain, and constant along every rectifiable curve in it; the
arbitrary domain is essential, and on a nowhere dense domain the time may vary
\cite{KoiralaPrep}. Vertical graphs are space-time graphs.

\begin{definition}
    A measure $\mu$ is \emph{regularly rectifiable of type $\Fam$} if $\mu\ll\HP$ and $\mu((\R^n\times\R)\setminus\bigcup_i\Gamma_i)=0$ for countably many graphs $\Gamma_i$ of the types allowed by $\Fam$, each with its own finite slope or norm.
\end{definition}

\subsection{The parabolic rectifiable Reifenberg theorem}

\begin{theorem}[Parabolic rectifiable Reifenberg]\label{thm:reifenberg}
For every $n\ge1$ there is $C=C(n)>0$ with the following property. Let $1\le k\le n+2$, let $\mu$ be a
locally Radon measure on $\R^n\times\R$, and let $M,\Gamma\ge0$ satisfy
\[
 \mu\bigl(P(0,1)\cap\{S^2_\mu>M\}\bigr)\le\Gamma .
\]
Then there are disjoint Borel sets $D_H,D_V\subset P(0,1)$ such that, writing $\mu_H=\mu\restr D_H$,
$\mu_V=\mu\restr D_V$, $\mu_+=\mu\restr(D_H\cup D_V)^c$ and $K_H=\supp\mu_H$, $K_V=\supp\mu_V$,
$K=K_H\cup K_V$:
\begin{enumerate}
\item $\mu=\mu_++\mu_H+\mu_V$ and $\mu_+(P(0,1))\le CM+\Gamma$;
\item $K_H$ and $K_V$ are compact subsets of $\overline{P(0,1)}$; $\mu_H=0$ unless $k\le n$, and
 $\mu_V=0$ unless $k\ge2$;
\item for every $\varepsilon>0$, $K_H$ is covered up to an $\HP$-null set by countably many horizontal
 graphs of slope at most $\varepsilon$, and $K_V$ by countably many regular vertical graphs of norm at
 most $\varepsilon$ with $G-c$ compactly supported for a constant $c$;
\item for all $q\in\R^n\times\R$ and $0<r\le R\le1$: $\HP(K\cap P(q,R))\le CR^k$; every $r$-separated subset of
 $K\cap P(q,R)$ has at most $C(R/r)^k$ points; $|\{z:\dist(z,K)<r\}|\le Cr^{n+2-k}$; and
 $\sum_ir_i^k\le CR^k$ for every finite or countable family of disjoint balls $P(z_i,r_i)$ with
 $z_i\in K\cap P(q,R)$ and $r_i\le R$.
\end{enumerate}
Moreover, the same constant works for a fixed given family: if
$\mu(P(0,1)\cap\{S^2_{\mu,\Fam}>M\})\le\Gamma$, the conclusion holds with $\mu_V=0$ when $\Fam=\Hor_k$ and
$\mu_H=0$ when $\Fam=\Ver_k$. For $k=n+2$ the conclusion holds with $\mu_+(P(0,1))=0$ and no hypothesis
on the beta numbers.
\end{theorem}

The Euclidean version of the rectifiable Reifenberg theorem for measures is due to Edelen, Naber and
Valtorta \cite[Theorem~1.1]{ENV2025}, in which a single family of $k$-planes gives one rectifiable part; see also \cite{ENV2019,NaberValtorta2017} and \cite[Theorem~8.2]{Naber2018}. In the parabolic setting, regular parabolic rectifiability of nodal and singular sets was proved in \cite{HallgrenKoiralaMa25} using quantitative uniqueness of tangent models. The competition between
the two types of models and the regularity in time of the vertical part are new here.

\begin{remark}\label{rem:converse}
\emph{The remainder is necessary.} Let $k<n+2$ and $\mu=c\,|\cdot|\restr P(0,1)$. Then
$\beta_\mu(z,r)^2\le Ccr^{n+2-k}$, so $S^2_\mu\le Cc$ and the hypothesis holds with $M=Cc$ and $\Gamma=0$.
Sets of finite $\HP$ measure are Lebesgue-null, so $\mu_H=\mu_V=0$ and all of $\mu$ is in the remainder.

\emph{There is no converse.} In $\R^2\times\R$, let $\ell_1=[-2,2]\times\{0\}\times\{0\}$ and
$\ell_2=[-2,2]\times\{h\}\times\{0\}$ with $0<h\le1/4$, let $k=1$ and $\mu=c(\HP[1]\restr\ell_1+\HP[1]\restr\ell_2)$.
Each $\ell_i$ is a horizontal graph of slope $0$, so for every $c>0$ the conclusion of
Theorem~\ref{thm:reifenberg} holds with $\mu_+(P(0,1))=0$. On the other hand, for $2h\le r\le1$ no single line
fits both segments in $P(z,r)$, and $\beta_\mu(z,r)^2\approx ch^2r^{-2}$; integrating over $r$ gives
$S^2_\mu\ge c_0c$ on $P(0,1)\cap\supp\mu$, with an absolute $c_0>0$. Hence the hypothesis
$\mu(P(0,1)\cap\{S^2_\mu>M\})\le\Gamma$ with $\Gamma<\mu(P(0,1))$ forces $M\ge c_0c$. Letting $c\to\infty$,
the decomposition gives no bound on $M$: beta numbers scale linearly in the mass, while the
decomposition does not see the mass. Under Ahlfors regularity, which fixes the scale of the mass,
the two sides become equivalent (Theorem~\ref{thm:quantitative}).
\end{remark}

\begin{corollary}\label{cor:density}
There is $C=C(n)>0$ such that for $1\le k\le n+2$ and a locally Radon $\mu$ with
$\mathcal E=\int_{P(0,1)}S^2_\mu\,d\mu<\infty$:
\begin{enumerate}
\item if $\Theta^k_*(\mu,\cdot)\le A$ $\mu$-a.e.\ on $P(0,1)$, then $\mu(P(0,1))\le C(A+\sqrt{\mathcal E})$;
\item if $\Theta^{k,*}(\mu,\cdot)\ge a>0$ $\mu$-a.e.\ on $P(0,1)$, then $\mu\restr P(0,1)$ is concentrated
 on the Borel set $F_a=P(0,1)\cap\{\Theta^{k,*}\ge a\}$, which for every $\varepsilon>0$ is covered by
 horizontal and regular vertical graphs with constants at most $\varepsilon$ and has
 $\HP(F_a)\le C(1+\sqrt{\mathcal E}/a)$;
\item under both assumptions $\mu\restr P(0,1)$ is regularly rectifiable.
\end{enumerate}
\end{corollary}

\subsection{Examples: sets that tilt in time}
The following examples show how
Theorem~\ref{thm:reifenberg} treats sets that look almost horizontal as discussed in \S \ref{sec:difficulties}.

\begin{figure}[ht]
\centering
\begin{tikzpicture}[>=Stealth,font=\scriptsize,
 declare function={chi(\s)=ifthenelse(abs(\s)<=1,1,ifthenelse(abs(\s)<2,cos(90*(abs(\s)-1))^2,0));}]
 \begin{scope}
  \draw[gray!60] (-1,-1) rectangle (1,1);
  \draw[->,gray] (-1.1,0)--(1.2,0) node[right] {$x$};
  \draw[->,gray] (0,-1.1)--(0,1.2) node[above] {$t$};
  \draw[very thick,setcol] (-1,-0.3)--(1,0.3);
  \draw[zoomcol,thick] (-0.1,-0.01) rectangle (0.1,0.01);
  \node at (0,-1.4) {(a) $P(0,1)$};
 \end{scope}
 \draw[zoomcol] (0.1,0.01)--(2.0,1);
 \draw[zoomcol] (0.1,-0.01)--(2.0,-1);
 \begin{scope}[xshift=3cm]
  \draw[zoomcol] (-1,-1) rectangle (1,1);
  \draw[->,gray] (-1.1,0)--(1.2,0) node[right] {$u$};
  \draw[->,gray] (0,-1.1)--(0,1.2) node[above] {$s$};
  \draw[very thick,setcol] (-1/3,-1)--(1/3,1);
  \node at (0,-1.4) {(b) $P(0,a/3)$ rescaled};
 \end{scope}
 \begin{scope}[xshift=8.4cm]
  \draw[->,gray] (-2.6,0)--(2.7,0) node[right] {$x$};
  \draw[->,gray] (0,-1.2)--(0,1.3) node[above] {$t$};
  \draw[very thick,setcol] (-2.4,-1.2)--(2.4,1.2);
  \foreach \p in {-1.8,-0.6,0.6,1.8} {
   \draw[densely dashed,modelcol] (\p,0.5*\p-0.45)--(\p,0.5*\p+0.45);
   \fill (\p,0.5*\p) circle (0.03);
  }
  \node at (0,-1.4) {(c) vertical models along the line};
 \end{scope}
\end{tikzpicture}
\caption{The slanted line $t=ax$ of Example~\ref{ex:slanted}. (a) In $P(0,1)$ it is close to the
horizontal axis (here $a=0.3$). (b) The small box $P(0,a/3)$ of (a), magnified by the parabolic
rescaling $\delta_{1/r}$, $r=a/3$: there the line is close to the vertical axis. (c) At the marked
points, the dashed lines are the vertical models $\{x_0\}\times\R$ that best approximate the line at
small scales (here $a=1/2$).}\label{fig:slanted}
\end{figure}
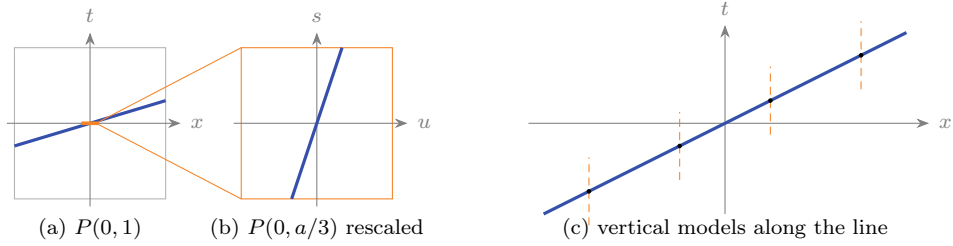

\begin{example}\label{ex:slanted}
Let $a>0$ and let $\mu_a$ be the image of the Lebesgue measure $\mathcal L^1\restr[-L,L]$ under
$x\mapsto(x,ax)$ in $\R\times\R$. For $k=1$ a direct computation gives
\[
 \beta_{\mu_a}(0,r)^2=\min\Bigl\{\frac ra,\frac ar\Bigr\},\qquad
 S^R_{\mu_a}(0)=2-\frac aR\quad(a\le R\le L).
\]
As $a\downarrow0$ the beta numbers tend to zero at every fixed scale, but $S^R_{\mu_a}(0)\to2$: the peak
$\beta^2=1$ moves to the scale $r=a$. Geometrically (Figure~\ref{fig:slanted}), the parabolic rescaling
$\delta_{1/r}(x,t)=(x/r,t/r^2)$ maps the line to $\{(u,(a/r)u)\}$. For $r\gg a$ this is within
parabolic distance $\sqrt{a/r}$ of the horizontal axis, while for $r\ll a$ it is within distance $r/a$
of the vertical axis. In particular the line has parabolic dimension two for every $a>0$.
\begin{itemize}
\item For $k=1$, the one-dimensional density of $\mu_a$ vanishes and every set of finite $\HP[1]$
 measure is $\mu_a$-null, so Theorem~\ref{thm:reifenberg} puts all of $\mu_a$ in the remainder $\mu_+$.
\item For $k=2$, the line is the vertical graph $x=t/a$. Fix a smooth cutoff $\chi$ with $\chi=1$ on
 $[-1,1]$ and $\chi=0$ outside $[-2,2]$. For $t_0\in\R$ and $r>0$, the function
 \[
  G_{t_0,r}(t)=\frac{t_0}a+\chi\Bigl(\frac{t-t_0}{r^2}\Bigr)\frac{t-t_0}a
 \]
 agrees with $t/a$ for $|t-t_0|\le r^2$, and its graph is a regular vertical graph with
 $\LipP G_{t_0,r}+\|\Dt G_{t_0,r}\|_{\BMOP}\le Cr/a$. Countably many such graphs cover the line, with
 arbitrarily small constants.
\end{itemize}
\end{example}

\begin{example}\label{ex:tilt}
In $\R\times\R$ let $C=\Gamma_-\cup\Sigma\cup\Gamma_+$ as in Figure~\ref{fig:tilt}, where
\[
 \Gamma_-=[-\tfrac34,-\tfrac14]\times\{\tfrac14\},\qquad
 \Sigma=\{(x,-x):|x|\le\tfrac14\},\qquad
 \Gamma_+=[\tfrac14,\tfrac34]\times\{-\tfrac14\}.
\]
On $\Sigma$ we have $\dP((x,-x),(y,-y))=|x-y|^{1/2}$, so $\Sigma$ has parabolic
dimension two and every subarc of $\Sigma$ has infinite $\HP[1]$ measure. The
segments $\Gamma_\pm$ have parabolic dimension one and are $\HP[2]$-null. Let $\ell$ be Euclidean
arclength on $C$. The decomposition of Theorem~\ref{thm:reifenberg} depends on $k$ and on the
measure, not only on the set; since $n=1$, there are no vertical models for $k=1$ and no horizontal
models for $k=2$. The table lists, in each case, a decomposition with the properties of
Theorem~\ref{thm:reifenberg}.

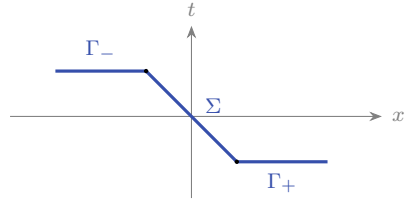
\begin{figure}[ht]
\centering
\begin{tikzpicture}[scale=2.4,>=Stealth,font=\scriptsize]
 \draw[->,gray] (-1,0)--(1.05,0) node[right] {$x$};
 \draw[->,gray] (0,-0.45)--(0,0.5) node[above] {$t$};
 \draw[very thick,setcol] (-0.75,0.25)--(-0.25,0.25) node[midway,above] {$\Gamma_-$};
 \draw[very thick,setcol] (-0.25,0.25)--(0.25,-0.25);
 \node[setcol] at (0.12,0.06) {$\Sigma$};
 \draw[very thick,setcol] (0.25,-0.25)--(0.75,-0.25) node[midway,below] {$\Gamma_+$};
 \foreach \x in {-0.25,0.25} \fill (\x,-\x) circle (0.012);
\end{tikzpicture}
\caption{The set $C$ of Example~\ref{ex:tilt}. The horizontal segments $\Gamma_\pm$ have parabolic
dimension one, while the ramp $\Sigma$ has parabolic dimension two: a parabolic ball of radius $r$
centered on $\Sigma$ contains a piece of $\Sigma$ of length $\approx r^2$.}\label{fig:tilt}
\end{figure}

\begin{center}
{\small
\begin{tabular}{@{}c l p{0.62\textwidth}@{}}
\toprule
$k$ & measure & decomposition\\ \midrule
$1$ & $\ell$ & $\mu_H=\ell\restr(\Gamma_-\cup\Gamma_+)$, $\mu_V=0$, $\mu_+=\ell\restr\Sigma$;
 $\supp\mu_H=\Gamma_-\cup\Gamma_+$. The ramp is always in the remainder, so $\ell(\Sigma)=\sqrt2/2$ must
 be at most $CM+\Gamma$.\\
$1$ & $\HP[1]\restr C$ & Not locally finite on $\Sigma$, so Theorem~\ref{thm:reifenberg} does not apply.\\
$2$ & $\ell$ & $\mu_H=0$, $\mu_V=\ell$, $\mu_+=0$; $\supp\mu_V=C$. The segments $\Gamma_\pm$ are
 $\HP[2]$-null exceptional sets of the vertical covers, so $\mu_V$ is not absolutely continuous
 with respect to $\HP[2]$.\\
$2$ & $\HP[2]\restr C$ & $\mu=\HP[2]\restr\Sigma$ lies on the regular vertical graph $x=G(t)$, with $G$
 smooth, compactly supported and $G(t)=-t$ for $|t|\le\frac14$. Canonical decomposition
 (Corollary~\ref{cor:canonical}): $\mu_V=\mu$, $\mu_H=\mu_+=0$, $\supp\mu=\Sigma$.\\
\bottomrule
\end{tabular}}
\end{center}
\end{example}

\section{Qualitative characterization}\label{sec:qualitative}

\begin{theorem}[Qualitative characterization]\label{thm:qualitative}
Let $1\le k\le n+2$, let $\Fam$ be a given family and $\mu$ a locally Radon measure on $\R^n\times\R$. The
following are equivalent:
\begin{enumerate}
\item $\mu$ is regularly rectifiable of type $\Fam$;
\item for some $R>0$, at $\mu$-a.e.\ point, $S^R_{\mu,\Fam}<\infty$, $\Theta^{k,*}(\mu,\cdot)>0$ and
 $\Theta^k_*(\mu,\cdot)<\infty$;
\item for some $R>0$, at $\mu$-a.e.\ point, $S^R_{\mu,\Fam}<\infty$ and $0<\Theta^{k,*}(\mu,\cdot)<\infty$;
\item $\mu\ll\HP$, $\mu$ is concentrated on a Borel set of $\sigma$-finite $\HP$ measure, and for some
 $R>0$, $S^R_{\mu,\Fam}<\infty$ $\mu$-a.e.
\end{enumerate}
In (2)--(4), ``for some $R>0$'' may be replaced by ``for every $R>0$''. If these hold, then
$0<\Theta^k_*\le\Theta^{k,*}<\infty$ $\mu$-a.e., and for every $\varepsilon>0$ the measure $\mu$ is
concentrated on countably many graphs of the types of $\Fam$ with slope or norm at most $\varepsilon$.
For $k=n+2$, $S^R_{\mu,\Fam}\equiv0$.
\end{theorem}

In $\R^d$, the characterization of rectifiability by the Jones square function under
$0<\Theta^{k,*}<\infty$ is due to Azzam--Tolsa (sufficiency, \cite[Theorem~1.1]{AzzamTolsa2015}) and
Tolsa (necessity, \cite[Theorem~1.1]{Tolsa2015}, \cite{Tolsa2019corr}); see also
\cite{Pajot1997,BadgerSchul2016,BadgerSchul2017} for related criteria. In the parabolic setting,
Bortz, Egert, Ferris and Saari characterized regular parabolic rectifiability of surface measure on a
codimension-one parabolic Lipschitz graph by a square function \cite[Theorem~1.3]{BEFS26}.
Theorem~\ref{thm:qualitative} applies to arbitrary measures, in all dimensions and for both types
of models. Regularity cannot be dropped: ordinary parabolic Lipschitz
graphs are not characterized by square functions (Example~\ref{ex:lacunary}; compare
\cite[Example~8.2]{Mat22}).

\begin{corollary}[Canonical decomposition]\label{cor:canonical}
Let $1\le k\le n+2$, and let $\mu$ be a locally Radon measure such that either $0<\Theta^{k,*}(\mu,\cdot)<\infty$
$\mu$-a.e., or $\mu\ll\HP$ and $\mu$ is concentrated on a Borel set of $\sigma$-finite $\HP$ measure. Then
$\mu$ decomposes uniquely as
\[
 \mu=\mu_H+\mu_V+\mu_+,
\]
where $\mu_H$ is regularly horizontally rectifiable, $\mu_V$ is regularly vertically rectifiable, and
$\mu_+$ vanishes on every horizontal graph and every regular vertical graph. The three measures are
locally Radon and mutually singular, and for every $R>0$,
\[
 \mu_H=\mu\restr\{S^R_{\mu,\Hor_k}<\infty\},\qquad \mu_V=\mu\restr\{S^R_{\mu,\Ver_k}<\infty\},\qquad
 \mu_+=\mu\restr\{S^R_{\mu,\Adm_k}=\infty\},
\]
with the convention that a set defined by an empty family is empty. In particular, $\mu_H$ and $\mu_V$
are covered by graphs of arbitrarily small slope or norm.
\end{corollary}

Writing $E_\Fam=\{S^R_{\mu,\Fam}<\infty\}$, the corollary contains the identities
$\mu(E_{\Hor_k}\cap E_{\Ver_k})=0$ and $\mu(E_{\Adm_k}\setminus(E_{\Hor_k}\cup E_{\Ver_k}))=0$. The second is not
a formal consequence of the definitions, since the best model in the mixed beta number may change type at every scale. It
reflects the separation of types: if $2\le k\le n$, a horizontal graph of slope $\varepsilon\le1$ and a
parabolic $L$-Lipschitz vertical graph with $L\varepsilon\le1/2$ meet in an $\HP$-null set, which has at
most one point when $k=2$. 

\begin{figure}[!htb]
\centering
\includegraphics[width=0.85\textwidth]{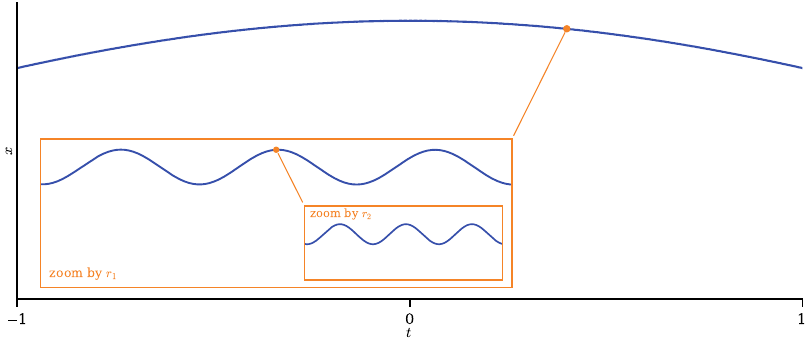}
\caption{The curve $F(t)=(f(t),t)$, $|t|\le1$, of Example~\ref{ex:lacunary}, drawn as $x=f(t)$ from
the formula truncated after three terms. The insets are the parabolic zooms: each reveals a new
oscillation, of amplitude $\varepsilon a_j$, and this continues at every scale.}\label{fig:lacunary}
\end{figure}

\begin{example}\label{ex:lacunary}
Let $a_j=(j+1)^{-1/2}$, $0<\varepsilon\le3/16$, and let $N_j$ be positive integers with
$N_{j+1}\ge16N_j$, $\pi N_j^{-1/2}\sum_{i<j}a_iN_i^{1/2}\le a_j/100$ and
$a_{j+1}N_{j+1}^{-1/2}\le a_jN_j^{-1/2}/1000$ (such $N_j$ exist by recursive choice). Put
\[
 f(t)=\varepsilon\sum_{j\ge0}a_jN_j^{-1/2}\cos(N_jt),\qquad F(t)=(f(t),t)\in\R\times\R,\qquad
 \nu=F_\#(\mathcal L^1\restr[-1,1]),
\]
where $\mathcal{L}^1$ is the Lebesgue measure, see Figure \ref{fig:lacunary}. Then, with $k=2$ and the mixed family, $F([-1,1])$ is compact with $0<\HP[2](F([-1,1]))<\infty$ and
is covered, for every $\varepsilon'>0$, by horizontal graphs and ordinary parabolic Lipschitz
vertical graphs of constant at most $\varepsilon'$. However $S^R_\nu(F(t))=\infty$ for $|t|<1$ and
every $R>0$, $\nu$ vanishes on every regular vertical graph, no countable family of horizontal and
regular vertical graphs covers $F([-1,1])$ up to an $\HP[2]$-null set, and $\nu$ is not regularly
rectifiable.
\end{example}

\section{Quantitative characterization}\label{sec:quantitative}

Let $A\ge1$. A measure $\mu$ is \emph{$A$-Ahlfors regular} if $\mu\ne0$ and, for all $z\in\supp\mu$,
\[
 \mu(P(z,r))\le Ar^k\quad(r>0),\qquad A^{-1}r^k\le\mu(P(z,r))\quad(0<r\le\diam\supp\mu).
\]
If $\supp\mu$ is unbounded, both bounds hold at all scales. If $\supp\mu$ is bounded, the lower bound is
only required up to its diameter. Since
$k\ge1$, the upper bound excludes atoms, so $\diam\supp\mu>0$. The
\emph{Carleson energy} of $\mu$ is
\[
 M_\beta(\mu)=\sup_{a\in\supp\mu,\ R>0}\frac1{R^k}\int_{P(a,R)}S^R_\mu(z)\,d\mu(z).
\]
For $L\ge1$, an \emph{$L$-reference} is one of the following: a \emph{horizontal} reference, the support
$\Gamma$ of $(F,T)_\#\mathcal L^k$ where $F:\R^k\to\R^n$ is $L$-bi-Lipschitz, $T:\R^k\to\R$ satisfies
$|T(u)-T(v)|\le L^2|u-v|^2$ on a nonempty compact set $D$, and $T$ is constant, equal to $T(a_Q)$ for a
point $a_Q\in D$ nearest to $Q$, on each cube $Q$ of a Whitney decomposition of $\R^k\setminus D$ into
dyadic cubes with $4\diam Q\le\dist(Q,D)\le10\diam Q$; a \emph{vertical} reference ($k<n+2$), the image $\Gamma$ of
$\Phi(u,t)=(F(u,t),t)$ where $\Phi:\R^{k-2}\times\R\to\R^n\times\R$ is $L$-bi-Lipschitz for $\dP$ and each
coordinate of $F$ has $\|\Dt F^\alpha\|_{\BMOP}\le L$; or, if $k=n+2$, $\Gamma=\R^n\times\R$.

\begin{definition}[Quantitative rectifiability]\label{def:qmr}
Let $K,L\ge1$. A \emph{$(K,L)$-covering} of $\mu$ on $P(a,R)$ consists of countably many
$L$-references $\Gamma_i$ and balls $P_i=P(a_i,R_i)\subset P(a,LR)$ with $a_i\in\supp\mu$ and
$P_i\cap\Gamma_i\neq\emptyset$, together with Borel sets $\mathcal B,\mathcal U_i$ partitioning
$(P(a,R)\cap\supp\mu)\times(0,R)$ up to sets of $\mu\otimes\frac{dr}r$-measure zero, such that
\[
 \mathcal U_i\subset\{(z,r):z\in P_i,\ L^{-1}\dist(z,\Gamma_i)\le r<R_i\},\qquad
 \iint_{\mathcal B}\frac{dr}{r}\,d\mu\le KR^k,
\]
\[
 \sum_iR_i^k\le KR^k,\qquad \sum_{a_i\in P(b,s),\,R_i\le s}R_i^k\le Ks^k\quad(b\in\R^n\times\R,\ s>0).
\]
The measure $\mu$ is \emph{$(K,L)$-quantitatively rectifiable} if $\mu\ll\HP$ and it has a
$(K,L)$-covering on every $P(a,R)$ with $a\in\supp\mu$, $R>0$.
\end{definition}

Thus every point--scale pair is either assigned to a reference visible at that scale, or to a bad
set of controlled logarithmic measure; this is a corona-type decomposition
\cite{DavidSemmes1991,DavidSemmes1993}.

\begin{theorem}[Quantitative characterization]\label{thm:quantitative}
Let $1\le k\le n+2$ and $A\ge1$. There are $L_0,C\ge1$ depending on $n,k,A$ such that every locally
Radon, $A$-Ahlfors regular measure $\mu$, with bounded or unbounded support, satisfies:
\begin{enumerate}
\item $M_\beta(\mu)<\infty$ if and only if $\mu$ is $(K,L_0)$-quantitatively rectifiable for some $K\ge1$;
\item if $M_\beta(\mu)<\infty$, then $\mu$ is $(C(1+M_\beta(\mu)),L_0)$-quantitatively rectifiable;
\item if $\mu$ is $(K,L_0)$-quantitatively rectifiable, then $M_\beta(\mu)\le CK$.
\end{enumerate}
Moreover, if $\supp\mu$ is unbounded, then for every $L\ge1$ there is $C(n,k,A,L)$ such that
$(K,L)$-quantitative rectifiability implies $M_\beta(\mu)\le C(n,k,A,L)K$.
\end{theorem}

In $\R^d$, following Jones \cite{Jones1990}, David and Semmes characterized uniform rectifiability by
Carleson bounds on beta numbers and by corona decompositions \cite{DavidSemmes1991,DavidSemmes1993}.
In codimension one, parabolic uniform rectifiability, defined by a Carleson condition on beta
numbers over vertical planes \cite{HLN03,HLN04}, was characterized by corona decompositions and big
pieces of regular graphs by Bortz, Hoffman, Hofmann, Luna-Garc\'ia and Nystr\"om
\cite[Theorem~1.2]{BHHLN23JGA}, and by big pieces of regular parabolic images by Bortz, Hyde and
Sharp \cite{BHS26}. In higher codimension, Hoffman and Jaye related vertical uniform rectifiability to the
boundedness of singular integrals \cite{HJ26}. Theorem~\ref{thm:quantitative} holds in every dimension, for both types of
models, with $L_0$ fixed before the size of $M_\beta$.

In both cases the Carleson energy and the coverings are taken over all radii $R>0$, and the
constants $L_0,C$ are the same. When the support is unbounded, the lower bound at all scales forces
$\mu(\R^n\times\R)=\infty$.

\begin{corollary}[Quantitative decomposition]\label{cor:approximants}
Let $\mu$ be as in Theorem~\ref{thm:quantitative}; if $\supp\mu$ is unbounded, assume $M_\beta(\mu)<\infty$.
Fix $K,L\ge1$ and a $(K,L)$-covering of $\mu$ on $P_0=P(a,R)$, $a\in\supp\mu$. Then for every $N\ge0$
there is a decomposition
\[
 \mu\restr P_0=\mu^N_H+\mu^N_V+e^N,\qquad e^N(P_0)\le2\exp\Bigl(-\frac{N}{(2^k+A^2)K}\Bigr)\mu(P_0),
\]
where $\mu^N_H$ is regularly horizontally rectifiable and $\mu^N_V$ is regularly vertically
rectifiable. As $N\to\infty$, $\mu^N_H$ and $\mu^N_V$ increase to the canonical
parts $\mu_H\restr P_0$ and $\mu_V\restr P_0$ of Corollary~\ref{cor:canonical}.
\end{corollary}

\begin{corollary}\label{cor:pieces}
A locally Radon measure $\mu$ is regularly rectifiable (mixed type) if and only if
$\mu=\sum_j\mu\restr E_j$ for disjoint Borel sets $E_j$ such that each nonzero $\mu\restr E_j$ is
$(2^k,\max\{2,\sqrt{1+n}\})$-quantitatively rectifiable. No Ahlfors regularity is assumed.
\end{corollary}

\section{Different exponents for beta numbers}\label{sec:betap}

For $p>0$ define the $L^p$ beta number $\beta_p(z,r)=\beta_{p,\mu,\Fam}(z,r)$ by
\[
 \beta_p(z,r)^p=\frac1{r^{k+p}}\inf_{L\in\Fam}\int_{P(z,r)}\dP(w,L)^p\,d\mu(w),
\]
and the two square functions
\[
 S^a_p(z)=\int_0^a\beta_p(z,r)^2\frac{dr}r,\qquad T^a_p(z)=\int_0^a\beta_p(z,r)^p\frac{dr}r ,
\]
so that $\beta_2=\beta$ and $S^a_2=T^a_2=S^a$.

\begin{corollary}\label{cor:table}
Let $1\le k<n+2$, $a>0$, and consider locally Radon $\mu$ with $0<\Theta^{k,*}(\mu,\cdot)<\infty$
$\mu$-a.e.\ and the mixed family. Finiteness of $S^a_p$ or $T^a_p$ almost everywhere is necessary or
sufficient for regular rectifiability as follows:
\begin{center}
\begin{tabular}{lccc}
\toprule
 & $0<p<2$ & $p=2$ & $p>2$\\ \midrule
$S^a_p<\infty$ a.e. & necessary, not sufficient & necessary and sufficient & sufficient, not necessary\\
$T^a_p<\infty$ a.e. & sufficient, not necessary & necessary and sufficient & necessary, not sufficient\\
\bottomrule
\end{tabular}
\end{center}
Each ``not'' means that the implication fails for some $(n,k)$: the counterexamples for $S_p$ have
$n=2$, $k=1$, and those for $T_p$ have $n=1$, $k=2$. For $k=n+2$ every such $\mu$ is regularly
rectifiable and all the energies vanish identically. Finiteness of $S^a_p$ and of $T^a_p$ does not
depend on $a$.
\end{corollary}

In $\R^d$, Tolsa showed that finiteness of the square function of $\beta_p$, $p\ne2$, does not
characterize rectifiability \cite[Theorems~1.1--1.2]{Tolsa2019beta}; for sufficient conditions with
$p\ne2$ under density assumptions see \cite{Pajot1997,BadgerSchul2016}. The parabolic counterexamples
for $S_p$ are planar ones placed in one time slice, which preserves the beta numbers exactly; those
for $T_p$ use lacunary graphs as in Example~\ref{ex:lacunary}.

\section{Lean verification}\label{sec:lean}

Every theorem and corollary of Sections~\ref{sec:decomposition}--\ref{sec:betap}, and
Example~\ref{ex:lacunary}, is proved in Lean~4 with Mathlib~\cite{deMouraUllrich21,Mathlib20}; the
file \texttt{docs/theorems.md} of the repository gives their Lean names. The formalized artifact is
the following.

\begin{center}
{\small
\begin{tabular}{@{}ll@{}}
\toprule
Repository & \url{https://github.com/robertkoirala/parabolic-rectifiability}\\
Release & tag \href{https://github.com/robertkoirala/parabolic-rectifiability/tree/v0.0.0}{\texttt{v0.0.0}}, commit \texttt{\releasecommit}\\
Lean / Mathlib & \texttt{v4.35.0-rc4} / commit \texttt{105b9b1ce126e77b467aac6f6ba36b510ac57368}\\
Size & 957 modules, about 176{,}000 lines, 14{,}534 declarations\\
Specification & \texttt{verification/Spec.lean}: 104 results and 18 auxiliary statements, Mathlib only\\
Axioms & \texttt{propext}, \texttt{Classical.choice}, \texttt{Quot.sound}\\
Checkers & Lean kernel; Comparator (\texttt{lake comparator}) with the nanoda kernel\\
\bottomrule
\end{tabular}}
\end{center}

\begin{sloppypar}\noindent Three checks run in continuous integration on Linux:
(i) \texttt{lake build} compiles every module, and Lean's kernel checks every proof, with no errors or
warnings; (ii) \texttt{lake env lean verification/AxiomAudit.lean} checks that every declaration of the
library uses only the three axioms above, so in particular there is no \texttt{sorry}; (iii) the file
\texttt{verification/Spec.lean}, which imports only Mathlib, restates the checked results and every
definition they use with proofs replaced by \texttt{sorry}, and Comparator checks that the library
proves exactly these statements, with only these axioms, replaying the proofs in Lean's kernel and
in the independent checker nanoda \cite{ComparatorRepo,NanodaLib}.\end{sloppypar}
To confirm that the theorems say what is claimed here, it suffices to read the statements and
definitions in \texttt{verification/Spec/}. The repository documents where the Lean encoding differs
from the conventions of a paper: for instance, energies take values in $[0,\infty]$, the infimum over
an empty family is $+\infty$, and graph covers are taken modulo $\HP$ outer measure.

\bibliographystyle{amsalpha}
{\footnotesize\bibliography{references}}

@article{AzzamTolsa2015,
  author = {Azzam, Jonas and Tolsa, Xavier},
  title = {{Characterization of $n$-rectifiability in terms of Jones' square function: Part II}},
  journal = {Geometric and Functional Analysis},
  volume = {25},
  number = {5},
  pages = {1371--1412},
  year = {2015},
  eprint = {1501.01572v2},
  archiveprefix = {arXiv},
  doi = {10.1007/s00039-015-0334-7}
}

@article{BadgerSchul2016,
  author = {Badger, Matthew and Schul, Raanan},
  title = {{Two sufficient conditions for rectifiable measures}},
  journal = {Proceedings of the American Mathematical Society},
  volume = {144},
  number = {6},
  pages = {2445--2454},
  year = {2016},
  eprint = {1412.8357v2},
  archiveprefix = {arXiv},
  doi = {10.1090/proc/12881}
}

@article{BadgerSchul2017,
  author = {Badger, Matthew and Schul, Raanan},
  title = {{Multiscale Analysis of 1-rectifiable Measures II: Characterizations}},
  journal = {Analysis and Geometry in Metric Spaces},
  volume = {5},
  number = {1},
  pages = {1--39},
  year = {2017},
  eprint = {1602.03823v3},
  archiveprefix = {arXiv},
  doi = {10.1515/agms-2017-0001}
}

@misc{Bamler20a,
  author = {Bamler, Richard H.},
  title = {{Structure Theory of Non-Collapsed Limits of {Ricci} Flows}},
  year = {2020},
  eprint = {2009.03243},
  archiveprefix = {arXiv},
  url = {https://arxiv.org/abs/2009.03243},
  note = {arXiv:2009.03243}
}

@misc{BEFS26,
  author = {Bortz, Simon and Egert, Moritz and Ferris, Sandra and Saari, Olli},
  title = {{The parabolic {D}ini-{$\beta$} condition and absolute continuity of surface and caloric measure}},
  year = {2026},
  eprint = {2608.23240},
  archiveprefix = {arXiv},
  primaryclass = {math.AP},
  note = {arXiv:2608.23240}
}

@article{Besicovitch1928,
  author = {Besicovitch, A. S.},
  title = {{On the fundamental geometrical properties of linearly measurable plane sets of points}},
  journal = {Mathematische Annalen},
  volume = {98},
  pages = {422--464},
  year = {1928},
  doi = {10.1007/bf01451603}
}

@article{BHHLN22AFM,
  author = {Bortz, Simon and Hoffman, John and Hofmann, Steve and Luna-Garcia, Jose Luis and Nystr{\"o}m, Kaj},
  title = {{On big pieces approximations of parabolic hypersurfaces}},
  journal = {Ann. Fenn. Math.},
  volume = {47},
  number = {1},
  year = {2022},
  pages = {533--571},
  doi = {10.54330/afm.115417}
}

@article{BHHLN22AIF,
  author = {Bortz, Simon and Hoffman, John and Hofmann, Steve and Luna-Garcia, Jose Luis and Nystr{\"o}m, Kaj},
  title = {{Coronizations and big pieces in metric spaces}},
  journal = {Ann. Inst. Fourier (Grenoble)},
  volume = {72},
  number = {5},
  year = {2022},
  pages = {2037--2078},
  doi = {10.5802/aif.3518}
}

@article{BHHLN23CME,
  author = {Bortz, Simon and Hoffman, John and Hofmann, Steve and Luna Garc{\'\i}a, Jos{\'e} Luis and Nystr{\"o}m, Kaj},
  title = {{Carleson measure estimates for caloric functions and parabolic uniformly rectifiable sets}},
  journal = {Anal. PDE},
  volume = {16},
  number = {4},
  year = {2023},
  pages = {1061--1088},
  doi = {10.2140/apde.2023.16.1061}
}

@article{BHHLN23JGA,
  author = {Bortz, S. and Hoffman, J. and Hofmann, S. and Luna-Garcia, J. L. and Nystr{\"o}m, K.},
  title = {{Corona decompositions for parabolic uniformly rectifiable sets}},
  journal = {J. Geom. Anal.},
  volume = {33},
  number = {3},
  year = {2023},
  pages = {Paper No. 96},
  doi = {10.1007/s12220-022-01176-8}
}

@misc{BHHW26,
  author = {Bortz, Simon and Hidalgo-Palencia, Pablo and Hofmann, Steve and Warta, James},
  title = {{Parabolic uniform rectifiability is characterized by {C}arleson measure estimates for bounded caloric functions}},
  year = {2026},
  eprint = {2609.15374},
  archiveprefix = {arXiv},
  primaryclass = {math.AP},
  note = {arXiv:2609.15374}
}

@article{BHMN25graph,
  author = {Bortz, Simon and Hofmann, Steve and Martell, Jos{\'e} Mar{\'\i}a and Nystr{\"o}m, Kaj},
  title = {{Solvability of the {$\mathrm{L}^{p}$} {D}irichlet problem for the heat equation is equivalent to parabolic uniform rectifiability in the case of a parabolic {L}ipschitz graph}},
  journal = {Invent. Math.},
  volume = {239},
  number = {1},
  year = {2025},
  pages = {165--217},
  doi = {10.1007/s00222-024-01300-1}
}

@misc{BHMN25pre,
  author = {Bortz, Simon and Hofmann, Steven and Martell, Jos{\'e} Mar{\'\i}a and Nystr{\"o}m, Kaj},
  title = {{Solvability of the {$L^p$} {D}irichlet problem for the heat equation implies parabolic uniform rectifiability}},
  year = {2025},
  eprint = {2510.22047},
  archiveprefix = {arXiv},
  primaryclass = {math.AP},
  note = {arXiv:2510.22047}
}

@misc{BHS26,
  author = {Bortz, Simon and Hyde, Matthew and Sharp, Mason},
  title = {{Big pieces of regular parabolic (bi-){L}ipschitz images is equivalent to parabolic uniform rectifiability}},
  year = {2026},
  eprint = {2608.23250},
  archiveprefix = {arXiv},
  primaryclass = {math.CA},
  note = {arXiv:2608.23250}
}

@article{CheegerHaslhoferNaber13,
  author = {Cheeger, Jeff and Haslhofer, Robert and Naber, Aaron},
  title = {{Quantitative Stratification and the Regularity of Mean Curvature Flow}},
  journal = {Geometric and Functional Analysis},
  volume = {23},
  number = {3},
  pages = {828--847},
  year = {2013},
  doi = {10.1007/s00039-013-0224-9}
}

@article{CheegerHaslhoferNaber14,
  author = {Cheeger, Jeff and Haslhofer, Robert and Naber, Aaron},
  title = {{Quantitative stratification and the regularity of harmonic map flow}},
  journal = {Calculus of Variations and Partial Differential Equations},
  volume = {53},
  number = {1-2},
  pages = {365--381},
  year = {2014},
  doi = {10.1007/s00526-014-0752-7}
}

@article{CheegerJiangNaber21,
  author = {Cheeger, Jeff and Jiang, Wenshuai and Naber, Aaron},
  title = {{Rectifiability of singular sets of noncollapsed limit spaces with Ricci curvature bounded below}},
  journal = {Annals of Mathematics},
  volume = {193},
  number = {2},
  pages = {407--538},
  year = {2021},
  doi = {10.4007/annals.2021.193.2.2}
}

@article{CheegerNaber12,
  author = {Cheeger, Jeff and Naber, Aaron},
  title = {{Lower bounds on Ricci curvature and quantitative behavior of singular sets}},
  journal = {Inventiones mathematicae},
  volume = {191},
  number = {2},
  pages = {321--339},
  year = {2012},
  doi = {10.1007/s00222-012-0394-3}
}

@article{CheegerNaber13,
  author = {Cheeger, Jeff and Naber, Aaron},
  title = {{Quantitative Stratification and the Regularity of Harmonic Maps and Minimal Currents}},
  journal = {Communications on Pure and Applied Mathematics},
  volume = {66},
  number = {6},
  pages = {965--990},
  year = {2013},
  doi = {10.1002/cpa.21446}
}

@article{CheegerNaber15,
  author = {Cheeger, Jeff and Naber, Aaron},
  title = {{Regularity of Einstein manifolds and the codimension 4 conjecture}},
  journal = {Annals of Mathematics},
  pages = {1093--1165},
  year = {2015},
  doi = {10.4007/annals.2015.182.3.5}
}

@article{CheegerNaberValtorta15,
  author = {Cheeger, Jeff and Naber, Aaron and Valtorta, Daniele},
  title = {{Critical Sets of Elliptic Equations}},
  journal = {Communications on Pure and Applied Mathematics},
  volume = {68},
  number = {2},
  pages = {173--209},
  year = {2014},
  doi = {10.1002/cpa.21518}
}

@article{ColdingMinicozzi16,
  author = {Colding, Tobias Holck and Minicozzi, William P.},
  title = {{The singular set of mean curvature flow with generic singularities}},
  journal = {Inventiones mathematicae},
  volume = {204},
  number = {2},
  pages = {443--471},
  year = {2015},
  doi = {10.1007/s00222-015-0617-5}
}

@misc{ComparatorRepo,
  author = {{Lean FRO and contributors}},
  title = {{{Comparator}: a trustworthy judge for {Lean} proofs}},
  howpublished = {GitHub repository, \url{https://github.com/leanprover/comparator}},
  year = {2026},
  url = {https://github.com/leanprover/comparator}
}

@book{DavidSemmes1991,
  author = {David, Guy and Semmes, Stephen},
  title = {{Singular integrals and rectifiable sets in $\mathbb{R}^n$. Au-del{\`a} des graphes lipschitziens}},
  series = {Ast{\'e}risque},
  number = {193},
  publisher = {Soci{\'e}t{\'e} Math{\'e}matique de France},
  year = {1991},
  url = {http://www.numdam.org/item/AST_1991__193__1_0/}
}

@book{DavidSemmes1993,
  author = {David, Guy and Semmes, Stephen},
  title = {{Analysis of and on Uniformly Rectifiable Sets}},
  series = {Mathematical Surveys and Monographs},
  volume = {38},
  publisher = {American Mathematical Society},
  year = {1993},
  doi = {10.1090/surv/038}
}

@article{DavidToro1999,
  author = {David, G. and Toro, T.},
  title = {{Reifenberg flat metric spaces, snowballs, and embeddings}},
  journal = {Mathematische Annalen},
  volume = {315},
  number = {4},
  pages = {641--710},
  year = {1999},
  doi = {10.1007/s002080050332}
}

@article{DavidToro2012,
  author = {David, Guy and Toro, Tatiana},
  title = {{Reifenberg parameterizations for sets with holes}},
  journal = {Memoirs of the American Mathematical Society},
  volume = {215},
  number = {1012},
  year = {2012},
  eprint = {0910.4869v1},
  archiveprefix = {arXiv},
  doi = {10.1090/s0065-9266-2011-00629-5}
}

@article{DeLellisMarcheseSpadaroValtorta18,
  author = {De Lellis, Camillo and Marchese, Andrea and Spadaro, Emanuele and Valtorta, Daniele},
  title = {{Rectifiability and upper Minkowski bounds for singularities of harmonic $Q$-valued maps}},
  journal = {Commentarii Mathematici Helvetici},
  volume = {93},
  number = {4},
  pages = {737--779},
  year = {2018},
  doi = {10.4171/cmh/449}
}

@incollection{deMouraUllrich21,
  author = {{de Moura}, Leonardo and Ullrich, Sebastian},
  title = {{The Lean 4 Theorem Prover and Programming Language}},
  booktitle = {Automated Deduction -- {CADE} 28},
  series = {Lecture Notes in Computer Science},
  volume = {12699},
  publisher = {Springer},
  pages = {625--635},
  year = {2021},
  doi = {10.1007/978-3-030-79876-5_37}
}

@article{Dorronsoro1985,
  author = {Dorronsoro, Jos{\'e} R.},
  title = {{A characterization of potential spaces}},
  journal = {Proceedings of the American Mathematical Society},
  volume = {95},
  number = {1},
  pages = {21--31},
  year = {1985},
  doi = {10.1090/s0002-9939-1985-0796440-3}
}

@article{EdelenEngelstein19,
  author = {Edelen, Nick and Engelstein, Max},
  title = {{Quantitative stratification for some free-boundary problems}},
  journal = {Transactions of the American Mathematical Society},
  volume = {371},
  number = {3},
  pages = {2043--2072},
  year = {2018},
  doi = {10.1090/tran/7401}
}

@article{ENV2019,
  author = {Edelen, Nick and Naber, Aaron and Valtorta, Daniele},
  title = {{Effective Reifenberg theorems in Hilbert and Banach spaces}},
  journal = {Mathematische Annalen},
  volume = {374},
  number = {3--4},
  pages = {1139--1218},
  year = {2019},
  eprint = {1806.01250v1},
  archiveprefix = {arXiv},
  doi = {10.1007/s00208-018-1770-0}
}

@article{ENV2025,
  author = {Edelen, Nick and Naber, Aaron and Valtorta, Daniele},
  title = {{Quantitative Reifenberg theorem for measures}},
  journal = {Mathematische Zeitschrift},
  volume = {310},
  number = {3},
  pages = {45},
  year = {2025},
  eprint = {1612.08052v3},
  archiveprefix = {arXiv},
  doi = {10.1007/s00209-025-03743-5},
  note = {Article number 45}
}

@misc{FangLi25a,
  author = {Fang, Hanbing and Li, Yu},
  title = {{Volume Estimates for the Singular Sets of Mean Curvature Flows}},
  year = {2025},
  eprint = {2504.09811},
  archiveprefix = {arXiv},
  url = {https://arxiv.org/abs/2504.09811},
  note = {arXiv:2504.09811}
}

@misc{FangLi25b,
  author = {Fang, Hanbing and Li, Yu},
  title = {{Singular Sets in Noncollapsed {Ricci} Flow Limit Spaces}},
  year = {2025},
  eprint = {2510.26317},
  archiveprefix = {arXiv},
  url = {https://arxiv.org/abs/2510.26317},
  note = {arXiv:2510.26317}
}

@article{Federer1947,
  author = {Federer, Herbert},
  title = {{The $(\varphi,k)$ rectifiable subsets of $n$-space}},
  journal = {Transactions of the American Mathematical Society},
  volume = {62},
  pages = {114--192},
  year = {1947},
  doi = {10.1090/s0002-9947-1947-0022594-3}
}

@article{FuWangWuZhang26,
  author = {Fu, Haotong and Wang, Wei and Wu, Ke and Zhang, Zhifei},
  title = {{Stratification and rectifiability of harmonic map flows via tangent measures}},
  journal = {Mathematische Annalen},
  volume = {395},
  number = {4},
  year = {2026},
  doi = {10.1007/s00208-026-03524-3}
}

@article{Gianniotis19,
  author = {Gianniotis, Panagiotis},
  title = {{Regularity theory for type I Ricci flows}},
  journal = {Calculus of Variations and Partial Differential Equations},
  volume = {58},
  number = {6},
  year = {2019},
  doi = {10.1007/s00526-019-1644-7}
}

@misc{HallgrenKoiralaMa25,
  author = {Hallgren, Max and Koirala, Robert and Ma, Zilu},
  title = {{Structure Theory of Parabolic Nodal and Singular Sets}},
  year = {2025},
  eprint = {2511.11570},
  archiveprefix = {arXiv},
  primaryclass = {math.AP},
  url = {https://arxiv.org/abs/2511.11570v3},
  note = {arXiv:2511.11570}
}

@misc{HJ26,
  author = {Hoffman, John and Jaye, Benjamin},
  title = {{On singular integrals and quantitative rectifiability in parabolic space and the {H}eisenberg group}},
  year = {2025},
  eprint = {2510.26934},
  archiveprefix = {arXiv},
  primaryclass = {math.AP},
  note = {arXiv:2510.26934}
}

@article{HL96,
  author = {Hofmann, Steve and Lewis, John L.},
  title = {{{$L^2$} solvability and representation by caloric layer potentials in time-varying domains}},
  journal = {Ann. of Math. (2)},
  volume = {144},
  number = {2},
  year = {1996},
  pages = {349--420},
  doi = {10.2307/2118595}
}

@article{HLN03,
  author = {Hofmann, Steve and Lewis, John L. and Nystr{\"o}m, Kaj},
  title = {{Existence of big pieces of graphs for parabolic problems}},
  journal = {Ann. Acad. Sci. Fenn. Math.},
  volume = {28},
  number = {2},
  year = {2003},
  pages = {355--384}
}

@article{HLN04,
  author = {Hofmann, Steve and Lewis, John L. and Nystr{\"o}m, Kaj},
  title = {{Caloric measure in parabolic flat domains}},
  journal = {Duke Math. J.},
  volume = {122},
  number = {2},
  year = {2004},
  pages = {281--346},
  doi = {10.1215/s0012-7094-04-12222-5}
}

@article{Hof97,
  author = {Hofmann, Steve},
  title = {{Parabolic singular integrals of Calder{\'o}n-type, rough operators, and caloric layer potentials}},
  journal = {Duke Math. J.},
  volume = {90},
  number = {2},
  year = {1997},
  pages = {209--259},
  doi = {10.1215/s0012-7094-97-09008-6}
}

@misc{HuangJiang23,
  author = {Huang, Yiqi and Jiang, Wenshuai},
  title = {{Volume Estimates for Singular Sets and Critical Sets of Elliptic Equations with {H{\"o}lder} Coefficients}},
  year = {2023},
  eprint = {2309.08089},
  archiveprefix = {arXiv},
  url = {https://arxiv.org/abs/2309.08089},
  note = {arXiv:2309.08089}
}

@misc{HuangJiang24,
  author = {Huang, Yiqi and Jiang, Wenshuai},
  title = {{The Nodal Sets of Solutions to Parabolic Equations}},
  year = {2024},
  eprint = {2406.05877},
  archiveprefix = {arXiv},
  url = {https://arxiv.org/abs/2406.05877},
  note = {arXiv:2406.05877}
}

@article{JiangNaber21,
  author = {Jiang, Wenshuai and Naber, Aaron},
  title = {{$L^2$ curvature bounds on manifolds with bounded Ricci curvature}},
  journal = {Annals of Mathematics},
  volume = {193},
  number = {1},
  pages = {107--222},
  year = {2021},
  doi = {10.4007/annals.2021.193.1.2}
}

@article{Jones1990,
  author = {Jones, Peter W.},
  title = {{Rectifiable sets and the Traveling Salesman Problem}},
  journal = {Inventiones Mathematicae},
  volume = {102},
  number = {1},
  pages = {1--15},
  year = {1990},
  doi = {10.1007/bf01233418}
}

@unpublished{KoiralaPrep,
  author = {Koirala, Robert},
  title = {{TBD}},
  year = {2026},
  note = {In preparation}
}

@article{KW80,
  author = {Kaufman, Robert and Wu, Jang-Mei},
  title = {{Singularity of parabolic measures}},
  journal = {Compositio Math.},
  volume = {40},
  number = {2},
  year = {1980},
  pages = {243--250},
  url = {https://www.numdam.org/item/CM_1980__40_2_243_0/}
}

@book{LM95,
  author = {Lewis, John L. and Murray, Margaret A. M.},
  title = {{The method of layer potentials for the heat equation in time-varying domains}},
  series = {Mem. Amer. Math. Soc.},
  volume = {114},
  number = {545},
  year = {1995},
  doi = {10.1090/memo/0545}
}

@article{LS88,
  author = {Lewis, John L. and Silver, Judy},
  title = {{Parabolic measure and the {D}irichlet problem for the heat equation in two dimensions}},
  journal = {Indiana Univ. Math. J.},
  volume = {37},
  number = {4},
  year = {1988},
  pages = {801--839},
  doi = {10.1512/iumj.1988.37.37039}
}

@misc{LW26,
  author = {Liu, Yu Tong and Workman, Myles},
  title = {{Parabolic rectifiability of the {B}rakke flow}},
  year = {2026},
  eprint = {2606.22441},
  archiveprefix = {arXiv},
  primaryclass = {math.DG},
  note = {arXiv:2606.22441}
}

@article{Marstrand1964,
  author = {Marstrand, John M.},
  title = {{The $(\varphi,s)$ regular subsets of $n$-space}},
  journal = {Transactions of the American Mathematical Society},
  volume = {113},
  pages = {369--392},
  year = {1964},
  doi = {10.1090/s0002-9947-1964-0166336-x}
}

@article{Mat22,
  author = {Mattila, Pertti},
  title = {{Parabolic rectifiability, tangent planes and tangent measures}},
  journal = {Ann. Fenn. Math.},
  volume = {47},
  number = {2},
  year = {2022},
  pages = {855--884},
  doi = {10.54330/afm.119821}
}

@incollection{Mathlib20,
  author = {{The mathlib Community}},
  title = {{The lean mathematical library}},
  booktitle = {Proceedings of the 9th ACM SIGPLAN International Conference on Certified Programs and Proofs},
  publisher = {ACM},
  pages = {367--381},
  year = {2020},
  doi = {10.1145/3372885.3373824}
}

@book{MatSurvey,
  author = {Mattila, Pertti},
  title = {{Rectifiability: A Survey}},
  publisher = {Cambridge University Press},
  year = {2023},
  doi = {10.1017/9781009288057}
}

@article{Mattila1975,
  author = {Mattila, Pertti},
  title = {{Hausdorff $m$ regular and rectifiable sets in $n$-space}},
  journal = {Transactions of the American Mathematical Society},
  volume = {205},
  pages = {263--274},
  year = {1975},
  doi = {10.1090/s0002-9947-1975-0357741-4}
}

@article{Miskiewicz2018,
  author = {Mi{\'s}kiewicz, Micha{\l}},
  title = {{Discrete Reifenberg-type theorem}},
  journal = {Annales Academiae Scientiarum Fennicae Mathematica},
  volume = {43},
  number = {1},
  pages = {3--19},
  year = {2018},
  eprint = {1612.02461v2},
  archiveprefix = {arXiv},
  doi = {10.5186/aasfm.2018.4301}
}

@article{MMP26,
  author = {Merlo, Andrea and Mourgoglou, Mihalis and Puliatti, Carmelo},
  title = {{On the density problem in the parabolic space}},
  journal = {Mem. Amer. Math. Soc.},
  volume = {320},
  number = {1630},
  year = {2026},
  doi = {10.1090/memo/1630}
}

@article{MPT22,
  author = {Mateu, Joan and Prat, Laura and Tolsa, Xavier},
  title = {{Removable singularities for {L}ipschitz caloric functions in time varying domains}},
  journal = {Rev. Mat. Iberoam.},
  volume = {38},
  number = {2},
  year = {2022},
  pages = {547--588},
  doi = {10.4171/rmi/1284}
}

@misc{Naber2018,
  author = {Naber, Aaron},
  title = {{Lecture notes on rectifiable {R}eifenberg for measures}},
  year = {2018},
  eprint = {1812.07564v1},
  archiveprefix = {arXiv},
  note = {arXiv:1812.07564v1}
}

@article{NaberValtorta17b,
  author = {Naber, Aaron and Valtorta, Daniele},
  title = {{Volume Estimates on the Critical Sets of Solutions to Elliptic PDEs}},
  journal = {Communications on Pure and Applied Mathematics},
  volume = {70},
  number = {10},
  pages = {1835--1897},
  year = {2017},
  doi = {10.1002/cpa.21708}
}

@article{NaberValtorta20,
  author = {Naber, Aaron and Valtorta, Daniele},
  title = {{The singular structure and regularity of stationary varifolds}},
  journal = {Journal of the European Mathematical Society},
  volume = {22},
  number = {10},
  pages = {3305--3382},
  year = {2020},
  doi = {10.4171/jems/987}
}

@article{NaberValtorta2017,
  author = {Naber, Aaron and Valtorta, Daniele},
  title = {{Rectifiable-Reifenberg and the regularity of stationary and minimizing harmonic maps}},
  journal = {Annals of Mathematics (2)},
  volume = {185},
  number = {1},
  pages = {131--227},
  year = {2017},
  eprint = {1504.02043v5},
  archiveprefix = {arXiv},
  doi = {10.4007/annals.2017.185.1.3}
}

@article {naber-valtorta-YM,
    AUTHOR = {Naber, Aaron and Valtorta, Daniele},
     TITLE = {Energy identity for stationary {Y}ang {M}ills},
   JOURNAL = {Invent. Math.},
  FJOURNAL = {Inventiones Mathematicae},
    VOLUME = {216},
      YEAR = {2019},
    NUMBER = {3},
     PAGES = {847--925},
      ISSN = {0020-9910,1432-1297},
   MRCLASS = {58E15 (53C07)},
  MRNUMBER = {3955711},
MRREVIEWER = {Graeme\ Wilkin},
       DOI = {10.1007/s00222-019-00854-9},
       URL = {https://doi-org.proxy.libraries.rutgers.edu/10.1007/s00222-019-00854-9},
}

@misc{naber-valtorta-2024energy,
  author = {Naber, Aaron and Valtorta, Daniele},
  title = {{Energy identity for stationary harmonic maps}},
  year = {2024},
  note = {arXiv:2401.02242}
}

@misc{NanodaLib,
  author = {{ammkrn and contributors}},
  title = {{nanoda\_lib: an external type checker for {Lean} 4}},
  howpublished = {GitHub repository, \url{https://github.com/ammkrn/nanoda_lib}},
  year = {2026},
  url = {https://github.com/ammkrn/nanoda_lib}
}

@article{NS17,
  author = {Nystr{\"o}m, Kaj and Str{\"o}mqvist, Martin},
  title = {{On the parabolic {L}ipschitz approximation of parabolic uniform rectifiable sets}},
  journal = {Rev. Mat. Iberoam.},
  volume = {33},
  number = {4},
  year = {2017},
  pages = {1397--1422},
  doi = {10.4171/rmi/976}
}

@article{Okikiolu1992,
  author = {Okikiolu, Kate},
  title = {{Characterization of subsets of rectifiable curves in $\mathbb{R}^n$}},
  journal = {Journal of the London Mathematical Society (2)},
  volume = {46},
  number = {2},
  pages = {336--348},
  year = {1992},
  doi = {10.1112/jlms/s2-46.2.336}
}

@article{Orp21,
  author = {Orponen, Tuomas},
  title = {{An integralgeometric approach to {D}orronsoro estimates}},
  journal = {Int. Math. Res. Not. IMRN},
  volume = {2021},
  number = {22},
  year = {2021},
  pages = {17170--17200},
  doi = {10.1093/imrn/rnz317}
}

@article{Pajot1997,
  author = {Pajot, Herv{\'e}},
  title = {{Conditions quantitatives de rectifiabilit{\'e}}},
  journal = {Bulletin de la Soci{\'e}t{\'e} Math{\'e}matique de France},
  volume = {125},
  number = {1},
  pages = {15--53},
  year = {1997},
  doi = {10.24033/bsmf.2298}
}

@article{Preiss1987,
  author = {Preiss, David},
  title = {{Geometry of measures in $\mathbb{R}^n$: distribution, rectifiability, and densities}},
  journal = {Annals of Mathematics (2)},
  volume = {125},
  number = {3},
  pages = {537--643},
  year = {1987},
  doi = {10.2307/1971410}
}

@article{Reifenberg1960,
  author = {Reifenberg, E. R.},
  title = {{Solution of the Plateau problem for m-dimensional surfaces of varying topological type}},
  journal = {Acta Mathematica},
  volume = {104},
  number = {1--2},
  pages = {1--92},
  year = {1960},
  doi = {10.1007/bf02547186}
}

@article{Schul2007,
  author = {Schul, Raanan},
  title = {{Subsets of rectifiable curves in Hilbert space---the analyst's TSP}},
  journal = {Journal d'Analyse Math{\'e}matique},
  volume = {103},
  pages = {331--375},
  year = {2007},
  eprint = {math/0602675v2},
  archiveprefix = {arXiv},
  doi = {10.1007/s11854-008-0011-y}
}

@article{Simon1993,
  author = {Simon, Leon},
  title = {{Cylindrical tangent cones and the singular set of minimal submanifolds}},
  journal = {Journal of Differential Geometry},
  volume = {38},
  number = {3},
  pages = {585--652},
  year = {1993},
  doi = {10.4310/jdg/1214454484}
}

@article{Simon1995,
  author = {Simon, Leon},
  title = {{Rectifiability of the singular set of energy minimizing maps}},
  journal = {Calculus of Variations and Partial Differential Equations},
  volume = {3},
  number = {1},
  pages = {1--65},
  year = {1995},
  doi = {10.1007/s005260050006}
}

@article{Tolsa2015,
  author = {Tolsa, Xavier},
  title = {{Characterization of $n$-rectifiability in terms of Jones' square function: Part I}},
  journal = {Calculus of Variations and Partial Differential Equations},
  volume = {54},
  number = {4},
  pages = {3643--3665},
  year = {2015},
  eprint = {1501.01569v3},
  archiveprefix = {arXiv},
  doi = {10.1007/s00526-015-0917-z}
}

@article{Tolsa2019beta,
  author = {Tolsa, Xavier},
  title = {{Rectifiability of measures and the $\beta_p$ coefficients}},
  journal = {Publicacions Matem{\`a}tiques},
  volume = {63},
  pages = {491--519},
  year = {2019},
  eprint = {1708.02304v3},
  archiveprefix = {arXiv},
  doi = {10.5565/publmat6321904}
}

@article{Tolsa2019corr,
  author = {Tolsa, Xavier},
  title = {{Correction to: Characterization of $n$-rectifiability in terms of Jones' square function: Part I}},
  journal = {Calculus of Variations and Partial Differential Equations},
  volume = {58},
  number = {1},
  pages = {35},
  year = {2019},
  doi = {10.1007/s00526-018-1478-8}
}

@incollection{Tolsa26ICM,
  author = {Tolsa, Xavier},
  title = {{Interactions between quantitative rectifiability, singular integrals, and boundary value problems for harmonic functions}},
  booktitle = {Proceedings of the International Congress of Mathematicians 2026, Vol. 2: Plenary Lectures},
  year = {2026},
  pages = {313--340},
  doi = {10.1137/25M1804054},
  eprint = {2607.16457},
  archiveprefix = {arXiv}
}

@article{Toro1995,
  author = {Toro, Tatiana},
  title = {{Geometric conditions and existence of bi-Lipschitz parameterizations}},
  journal = {Duke Mathematical Journal},
  volume = {77},
  number = {1},
  pages = {193--227},
  year = {1995},
  doi = {10.1215/s0012-7094-95-07708-4}
}

\end{document}